\documentclass[11pt,leqno]{article}
\usepackage{amssymb}
\usepackage{latexsym}
\usepackage{amsfonts, amsmath}

\newcommand{\demo}{ \noindent {\it   Proof. }}
\newcommand{\qed}{$\Box$\bigskip}

\newcommand{\eqg}{=_{_G}}
\newcommand{\la}[1]{|#1|_{_A}}
\newcommand{\lf}[1]{|#1|_{_F}}
\newcommand{\lgr}[1]{|#1|_{_G}}
\newcommand{\lh}[1]{|#1|_{_H}}
\newcommand{\lz}[1]{|#1|_{_{\mathbb{Z}}}}
\newcommand{\lzr}[1]{|#1|_{_{\mathbb{Z}^r}}}

\newcommand{\ar}{\mathcal{A}}
\newcommand{\n}{\mathcal{N}}

\newtheorem{thm}{Theorem}[section]
\newtheorem{prop}[thm]{Proposition}
\newtheorem{lem}[thm]{Lemma}
\newtheorem{cor}[thm]{Corollary}

\newtheorem{ass}[thm]{Assumption}

\newcommand{\area}{\operatorname{area}}
\newcommand{\var}{\operatorname{var}}

\begin{document}

\title{The mean Dehn functions of abelian groups}
\author{O. Bogopolski$^{(1)}$ \hspace{1.5cm}  and  \hspace{1.5cm} E. Ventura$^{(2)}$ \vspace{.6cm} \\
\begin{tabular}{lcl}
$^{(1)}$Institute of Mathematics, & & $^{(2)}$Dept.\ Mat.\ Apl.\ III, \\
Koptjuga 4, Novosibirsk, & & Univ.\ Pol.\ Catalunya, \\
630090, Russia & & Barcelona, Catalunya (Spain) \\
\phantom{aaaaaaaa} and  & & \phantom{aaaaaaaa} and \\
Universit\"{a}t Dortmund & \phantom{a and a} & Centre de Recerca Matem\`atica \\
Fachbereich Mathematik &  & Barcelona, Catalunya (Spain) \\
Lehrstuhl VI (Algebra) &  &  \\
Vogelpothsweg 87 & & \\
D-44221 Dortmund, Germany & & \\
& & \\
& & \\
groups@math.nsc.ru & & enric.ventura@upc.edu
\end{tabular}
}

\date{\today}
\maketitle

\begin{abstract}
While Dehn functions, $D(n)$, of finitely presented groups are very well studied in the literature, mean Dehn functions
are much less considered. M. Gromov introduced the notion of mean Dehn function of a group, $D_{mean}(n)$, suggesting
that in many cases it should grow much more slowly than the Dehn function itself. Using only elementary counting
methods, this paper presents some computations pointing into this direction. Particularizing them to the case of any
finite presentation of a finitely generated abelian group (for which it is well known that $D(n)\sim n^2$ except in the
1-dimensional case), we show that the three variations $D_{osmean}(n)$, $D_{smean}(n)$ and $D_{mean}(n)$ all are
bounded above by $Kn(\ln n)^2$, where the constant $K$ depends only on the presentation (and the geodesic combing)
chosen. This improves an earlier bound given by Kukina and Roman'kov.
\end{abstract}

\section{Introduction}

For all the paper, let $A=\{ a_1,\, \ldots ,\, a_r\}$ be an alphabet with $r$ letters and let $A^{\ast}$ be the free
monoid on $A\cup A^{-1}$.

Let also $G$ be an $r$-generated finitely presented group, and choose a finite presentation $G=\langle A\,|\,
R\rangle$, with $A$ as set of generators. We have the natural epimorphisms $A^{\ast}\twoheadrightarrow
F\twoheadrightarrow G$, where $F$ is the free group on $A$. Whenever clear from the context, we shall use the same
notation for referring to a formal word $w$ in $A^{\ast}$, and to its images in $F$ and $G$. When necessary, we shall
use $=_{_A}$, $=_{_F}$ and $\eqg$ to denote equality in these three algebraic structures.

Let $w\in A^{\ast}$ be a (possibly non-reduced) word. We shall denote by $\la{w}$, $\lf{w}$ and $\lgr{w}$ the metric
lengths of $w$ in $A^{\ast}$, $F$ and $G$, respectively. In other words, $\la{w}$ equals the number of letters in $w$,
$\lf{w}$ means the number of letters in $w$ after free reduction, and $\lgr{w}$ equals the number of letters in the
shortest word $w'\in A^{\ast}$ such that $w\eqg w'$. To avoid possible confusions with lengths, we shall write the
cardinal of a set $S$ as $\sharp S$.

Clearly, if $H$ is a quotient of $G$, say $A^{\ast}\twoheadrightarrow F\twoheadrightarrow G\twoheadrightarrow H$, then
$\la{w}\geqslant \lf{w}\geqslant \lgr{w}\geqslant \lh{w}$. For example, taking $A=\{a\}$, $F=\langle a\rangle \simeq
\mathbb{Z}$, $G=\langle a\, |\, a^{10} \rangle \simeq \mathbb{Z}/10\mathbb{Z}$, $H=\langle a\, |\, a^5 \rangle \simeq
\mathbb{Z}/5\mathbb{Z}$ and $w=aaa^{-1}aaa$, we have $\la{w}=6$, $\lf{w}=4$, $\lgr{w}=4$ and $|w|_{_H}=1$.

Let $\Gamma(G)$ denote the Cayley graph of $G$ with respect to $A$, and let $e$ be the vertex corresponding to the
trivial element. There is a natural bijection, $w\longleftrightarrow \gamma_w$, between (possibly non-reduced) words in
$A^{\ast}$ and paths in $\Gamma(G)$ starting at $e$ (and possibly with backtrackings). In the future, we will not
distinguish between $w$ and $\gamma_w$, usually using $w$ to denote the corresponding path as well (if there is no risk
of confusion). Clearly, the length of $\gamma_w$ is $\la{w}$, the length of $\gamma_w$ after reducing all possible
backtrackings is $\lf{w}$, and the distance in $\Gamma(G)$ from $e$ to $\tau \gamma_w$ (the terminal point of
$\gamma_w$) is $\lgr{w}$. Any path in $\Gamma(G)$ of the minimal possible length from $e$ to $\tau \gamma_w$ is called
a \emph{geodesic for $w\in G$} and, in fact, it represents a word $w'\in A^{\ast}$ of the shortest possible $A$-length
such that $w\eqg w'$. Of course, geodesics are not unique, in general.

Let $w\in A^{\ast}$. Clearly, $w\eqg 1$ if and only if $\gamma_w$ is closed. In this case, $w\in F$ belongs to the
kernel of the projection $F\twoheadrightarrow G$ and so, it can be expressed as
 $$
w=\prod_{i=1}^m f_i^{-1}r_i^{\epsilon_i}f_i,
 $$
where $f_i\in F$, $r_i\in R$, and $\epsilon_i =\pm 1$. The minimal such $m$ is called the {\it area} of $w$, denoted
$\area(w)$. The motivation for this name is obviously of geometric nature. For every vertex $v\in \Gamma(G)$ and every
relator $r_i$, there is a closed path at $v$ which labels $r_i$. For every such path $p$, let us add a 2-cell to
$\Gamma(G)$ with boundary $p$. In the resulting 2-complex, the area of $w$ is the minimal number of 2-cells needed to
fill a disc with boundary $w$.

Note that if $w, w'\in A^{\ast}$ reduce to the same element in $F$ which maps to the identity element in $G$, then
$\area(w)=\area(w')$. It is clear from the definition that, for $w,w'\in A^{\ast}$ with $w\eqg w'\eqg 1$, we have
$\area(ww')\leqslant \area(w)+\area(w')$. Also, $\area(w^{-1})=\area(w)$ and $\area(vwv^{-1})=\area(w)$ for every $v\in
A^{\ast}$.

The way those areas grow when considering longer and longer words in the group $G$, is measured by the so-called Dehn
function associated to the prefixed presentation for $G$. To give the precise definition, we need the following
notation. For every positive integer $n$ define the sets
 $$
B_G(n)=\{ w\in A^{\ast} \mid w\eqg 1 ,\, \, \la{w}\leqslant n \}.
 $$
and
 $$
S_G(n)=\{ w\in A^{\ast} \mid w\eqg 1 ,\,\, \la{w}=n \} =B_G (n)\setminus B_G (n-1).
 $$
By convention, let us write $B_G (0)=S_G (0)=\{1\}$. The notation $B_G(n)$ and $S_G(n)$ reflects the idea of
\textit{balls} and \textit{spheres}, respectively. However, note that these sets are not real balls or spheres in the
metric of $G$, but sets of closed paths at $e$ with possible backtrackings, and  with bounded or given $A$-length.

Note that, if $H$ is a quotient of $G$ then $B_G (n) \subseteq B_H (n)$ and $S_G (n)\subseteq S_H (n)$. So, the bigger
sets correspond to the trivial group (in this case we delete the subindex to avoid confusions). This way,
 $$
B (n)=\{ w\in A^{\ast} \mid \la{w}\leqslant n \}
 $$
and
 $$
S (n)=\{ w\in A^{\ast} \mid \la{w}=n \}
 $$
are the real ball and the real sphere in the monoid $A^*$, respectively. Furthermore, it is easy to see that $\sharp
S_G (n) \leqslant \sharp S (n)=(2r)^n$ and $\sharp B_G (n) \leqslant \sharp B (n)=(2r)^0+(2r)^1+\cdots +(2r)^n
=\frac{(2r)^{n+1}-1}{2r-1}$.

Now, the {\it Dehn function} of the finite presentation $G=\langle A\,|\, R\rangle$ is the function
$D:\mathbb{N}\rightarrow \mathbb{N}$ defined by
 $$
D(n)= \underset{w\in B_G(n)}{\max} \{\area(w) \}.
 $$
It measures the biggest area of those words in the ball of a given radius. In principle, this function depends on the
presentation but it is well-know that, changing to another presentation of the same group, $D(n)$ remains the same up
to multiplicative and additive constants, both in the argument and in the range. In particular, the asymptotic behavior
of $D(n)$ only depends on $G$.

There are a lot of papers in the literature investigating Dehn functions of groups (specially because of its relation
with the word problem of the group). For example, it is well known that every word-hyperbolic group has a linear Dehn
function, and that automatic groups have Dehn function at most quadratic (see~\cite{E} for a general exposition). Also,
a relevant theorem attributed to Gromov states that every subquadratic Dehn function is in fact linear (see~\cite{O}
for a detailed proof), thus existing a gap between $n$ and $n^2$ on the asymptotic behavior of Dehn functions of
finitely generated groups. A consequence of these results is that non-cyclic finitely generated free abelian groups (as
automatic but non word-hyperbolic groups) have precisely quadratic Dehn function, i.e., $C_1 n^2 \leqslant
D(n)\leqslant C_2 n^2$ for appropriate constants $C_1, C_2 >0$.
%
%Here, we note only the following result:
%
%\begin{thm}[Gersten, Holt, Riley~\cite{GHR}] The Dehn function of any finitely generated nilpotent group of class $c$ is
%$O(n^{c+1})$.
%\end{thm}

\medskip

In the literature, there are interesting variations of the concept of Dehn function, which are still not deeply
investigated. In this paper, we are concerned to \emph{mean Dehn functions}, first introduced by M. Gromov in~\cite{G}.

The \emph{mean Dehn function} of the presentation $\langle A\,|\, R\rangle$ for $G$, denoted $D_{\text{mean}}$, is the
mapping $D_{\text{mean}} \colon \mathbb{N}\to \mathbb{Q}$ defined by
 $$
D_{\text{mean}}(n)=\frac{\underset{w\in B_G (n)}{\sum} \area(w)}{\sharp B_G (n)}
 $$
(note that the denominator is never $0$ since the empty word always belongs to $B_G (n)$).

Similarly, the {\it spherical mean Dehn function}, denoted $D_{\text{smean}}$, is defined as
 $$
D_{\text{smean}}(n)=\frac{\underset{w\in S_G (n)}{\sum} \area(w)}{\sharp S_G (n)},
 $$
where we understand $D_{\text{smean}}(n)=0$ if the sphere $S_G (n)$ is empty.

Since areas of words (and also balls and spheres) do depend on the chosen presentation for $G$, the functions
$D_{\text{mean}}$ and $D_{\text{smean}}$ also depend on that presentation. Contrasting with the situation for the
classical Dehn function, it is still not known in general whether the asymptotic behavior of these averaged versions is
also invariant under changing the presentation.

As we said, these averaged Dehn functions are still very poorly considered in the literature. One of the few existing
results is due to E. G. Kukina and V. A. Roman'kov~\cite{KR} who proved that, for finitely generated free abelian
groups,
 $$
\lim_{n\to \infty} \frac{D_{\text{mean}}(n)}{n^{7/4}} =0.
 $$
This is considerably improved in the present paper, where we give the following much better asymptotic bound:

\begin{thm}\label{main-ab}
The mean Dehn function of a finitely generated abelian group $G$ satisfies $D_{\text{mean}}(n)=O\bigl(n(\ln n)^2\bigr)$
(with the constant depending only on the chosen finite presentation for $G$). The same assertion is valid for the
spherical mean Dehn function of $G$.
\end{thm}

Here, as in the rest of the paper, we make use of the ``$O$" notation for comparing the growth of pairs of functions.
Given two functions $f,g\colon \mathbb{N}\to \mathbb{R}^+$ defined on the set of natural numbers and having positive
values, one writes $f(n)=O(g(n))$ when there exists a constant $K$ (independent on $n$) such that $f(n)\leqslant Kg(n)$
for every $n\geqslant 1$. Note that, by changing $K$ to $\max \{K,\, f(1)/g(1),\, \ldots ,\, f(n_0)/g(n_0)\}$, this is
the same as having the inequality for big enough $n$, say $n>n_0$ (we shall refer to this by writing $n\gg 0$). This
notation is useful when one is mostly interested on the existence of such constant, more than on its actual value (this
is the case in the present paper; however, following the details in our arguments, one can always extract from them a
concrete value for the corresponding constant in all our statements involving ``$O$").

Back to Dehn functions, it is worth remarking that M. Sapir recently introduced another interesting variation of the
concept of Dehn function, namely his \emph{random Dehn function}. It uses the notion of area of a word $w\in A^{\ast}$
(not necessarily equal to 1 in $G$) with respect to a given geodesic combing in $\Gamma(G)$ (see the beginning of
section~\ref{s-combings}, below). Having chosen a geodesic combing in $\Gamma(G)$, say that $f\colon \mathbb{N}\to
\mathbb{N}$ is a \emph{random isoperimetric function} for $G$ if
 $$
\frac{\sharp \{w\in A^{\ast}\,|\, \la{w}\leqslant n,\,\, {\text {\rm area}}(w)\geqslant f(n)\}}{\sharp \{ w\in A^{\ast}
\mid \la{w}\leqslant n\}}\rightarrow 0,
 $$
for $n\rightarrow \infty$. Then, the \emph{random Dehn function} for $G$ is the smallest random isoperimetric function
(which, a priori, depends on the presentation of $G$ and on the chosen combing). M. Sapir claimed (private
communication) that, for any finite presentation of an abelian group $G$, and for any geodesic combing in $\Gamma(G)$,
there exists a constant $K$ such that the random Dehn function of $G$ is dominated by $n\mapsto Kn\ln n$. It would be
interesting to investigate the possible relationships between mean and random Dehn functions.

To conclude this introduction, let us avoid possible notational confusions by saying that, all over the paper, we use
the term ``$\ln$" meaning neperian logarithm (i.e. $\exp(\ln n)=n$). Just for technical reasons ($\ln 1=0$ and we will
need to work with functions $f\colon \mathbb{N}\to \mathbb{R}^+$ taking strictly positive values) the set $\mathbb{N}$
will be taken to be all natural numbers except 1. Also, for every real number $x$, we shall denote by $\lfloor
x\rfloor$ its integral part (i.e. the biggest integer which is less than or equal to $x$), and $\lceil x\rceil =\lfloor
x\rfloor +1$. So, $\lfloor x\rfloor \leqslant x<\lceil x\rceil$. Note that, for a positive integer $n$, $n>x$ is
equivalent to $n\geqslant \lceil x\rceil$; and $n\leqslant x$ is equivalent to $n\leqslant \lfloor x\rfloor$. Also, for
every integer $n>0$, $\lfloor \frac n2 \rfloor + \lceil \frac n2 \rceil =n$.

\medskip

The paper is organized as follows. In Section~\ref{s-combings} we introduce the notion of open mean Dehn function and
give a general upper bound for it, assuming that the presentation satisfies some technical assumptions. We also give
some indications on how to convert this bound into a bound for the spherical mean and the mean Dehn functions. In
Section~\ref{s-counting} we concentrate on finitely generated abelian groups, making the necessary countings there to
ensure that every finite presentation of such a group satisfies the assumptions required in the previous section.
Finally, in Section~\ref{s-ab-case} we deduce explicit upper bounds for the mean and the spherical mean Dehn functions
of any finite presentation of an abelian group. It is interesting to remark that the techniques developed in
Section~\ref{s-combings} can probably be applied to other groups as well. As soon as one can find two functions
satisfying assumption~\ref{assumption} for his favorite group presentation, an upper bound for the open spherical mean
Dehn function of that presentation will follow easily. With some more computations, one can also hope to obtain an
upper bound for the mean Dehn function of such presentation.

\medskip

We have to mention that, during the long process of publication of the present paper, another preprint appeared with
similar results. Totally independently from us, R. Young~\cite{Y} considers finitely generated nilpotent groups and
proves several results about what he calls their \textit{averaged Dehn function}. His results imply that, for the
finitely generated abelian case, this function is $O(n\ln n)$. However, a rather technical but quite important detail
needs to be highlighted when comparing both papers (i.e. when comparing the definitions of averaged and mean Dehn
functions). In~\cite{Y}, the author considers what he calls \textit{lazy words}, which are elements of the free monoid
on $A\cup A^{-1}\cup \{ e\}$ i.e., formal sequences of the form $a_1 \cdots a_n$ with $a_i \in A\cup A^{-1}\cup \{
e\}$. Because of the possibility of using the symbol $e$ (which represents the trivial element in $G$), a lazy word of
length $n$ corresponds to a (non-necessarily reduced) word of length less than or equal to $n$, in our terminology. But
when counting them (and averaging their areas) there is a significant difference. The total number of lazy words of
length $n$ is $(2r+1)^n$, while the total number of our words of length less than or equal to $n$ is
$\frac{(2r)^{n+1}-1}{2r-1}$, asymptotically like $(2r)^n \ll (2r+1)^n$. The difference is due to the fact that every
word $w$ of length $m<n$ appears many times counted as a lazy word, precisely as many as ways there are of expanding
$w$ to a sequence of $n$ symbols by adding $n-m$ $``e"$'s between the existing ones. And all these different
representations of the same element of $G$, of course have the same area. So, for sure, this effect introduces an
artificial distortion when estimating the corresponding areas. When averaging the areas of lazy words of length $n$ (as
is done in~\cite{Y}) one is counting shorter words with bigger multiplicity (the maximum distortion appears around
words of length $n/2$). And, of course, those shorter words have smaller area in average. So, this distortion in the
counting contributes to artificially decrease the global average of areas. It is very difficult to make a quantitative
estimation of this effect, but we believe it can very well be the reason of the difference between the bound $O(n\ln
n)$ obtained in~\cite{Y}, and the bound $O(n(\ln n)^2)$ obtained here.

Beyond this discussion, there is the question of which is the good (...or the most appropriate, or the best ...) notion
of mean Dehn function from the group theory point of view. In other words, which is the exact set that must be
considered to average the areas over it? The appendix at the end of this paper pretends to contribute to this
discussion.

\section{Combings in groups and the open mean Dehn function}\label{s-combings}

For technical reasons, we will need an extension of the concept of area to arbitrary paths in $\Gamma(G)$ (not just
those which are closed at $e$, i.e. words in $A^{\ast}$ mapping to 1 in $G$). Accordingly, we shall introduce the
notion of open mean Dehn function averaging over all those words.

A \emph{combing} in $\Gamma(G)$ is a set $T$ consisting of exactly one path from $e$ to every vertex $v\in \Gamma(G)$,
denoted $T[e,\, v]$ or simply $T[v]$, and such that $T[e]$ is the trivial path. By translation, such a set also
determines a (unique) path between every given pair of vertices in $\Gamma(G)$, namely $T[u,\, v]=uT[e,\, u^{-1}v]$. A
combing $T$ is said to be \textit{geodesic} if $T[v]$ (and so, $T[u,v]$) is a geodesic path, for every pair of vertices
$u,v$. Using a combing $T$, any path $\gamma$ in $\Gamma(G)$ can be \textit{closed up} by returning back to its initial
vertex through the combing. That is, defining $\widetilde{\gamma }=T[\iota \gamma,\, \tau \gamma]$, we have that
$\gamma \widetilde{\gamma}^{-1}$ is a closed path at $\iota \gamma$. Note that if $T$ is geodesic then
$\la{\widetilde{\gamma}}\leqslant \la{\gamma}$.

Standard examples of combings are the \textit{tree combings}, i.e. those determined by a maximal tree $T$ in $\Gamma
(G)$. In this case, $T[v]$ is the unique reduced path from $e$ to $v$ in $T$. For example, $\Gamma(\mathbb{Z}^2)$ (with
the standard presentation) is the two dimensional integral lattice; and the maximal tree given by the $X$-axis plus all
the vertical lines, determines the geodesic combing of $G=\mathbb{Z}^2$ where $T[(r,s)]$ is the path that goes first
$r$ steps to the right and then $t$ steps up. Note that, for these tree combings, usually $T[wu,\, wv]=wT[u,\, v]$
\textit{is not} the path determined by the tree from $wu$ to $wv$.

With the help of combings, we can define the area of an arbitrary path $\gamma$ in $\Gamma(G)$ (not-necessarily
reduced, neither closed, neither even starting at $e$). If $\gamma$ is closed at $e$ we already know the meaning of
$\area(\gamma )$. If $\gamma$ is closed at a vertex $u=\iota \gamma =\tau \gamma \neq e$ we define the \textit{area} of
$\gamma$ by first translating $\gamma$ to $e$ (i.e. reading the same word $\gamma$ but from the vertex $e$) or,
equivalently, going first to (and then coming back from) $u$ through an arbitrary path (which makes no difference at
the level of the area because it is conjugacy invariant):
 $$
\area(\gamma)=\area(T[e,u]\gamma T[e,u]^{-1})
 $$
(caution! $T[e,u]^{-1}\neq T[u,e]=uT[e,u^{-1}]$ in general). Finally, suppose $\gamma$ is an arbitrary path in
$\Gamma(G)$ (with $u=\iota \gamma$ and $v=\tau \gamma$ not necessarily equal, neither equal to $e$). The \textit{area}
of $\gamma$ is defined by first closing it through the combing:
 $$
\area(\gamma)=\area(\gamma \widetilde{\gamma}^{-1}).
 $$
Since, by definition $T[u,v]=uT[e,u^{-1}v]$, closing up $\gamma$ and translating the result to $e$ reads the same as
translating first $\gamma$ to $e$ and then closing it up.

%Now, let $w$ be an arbitrary word in $A^{\ast}$ thought of as a (non-necessarily closed neither reduced) path in
%$\Gamma(G)$ starting at $e$. We define the \emph{area} of $w$ (with respect to a chosen combing) as the area of the
%closed path $w\widetilde{w}^{-1}$ at $e$, which is now a genuine area because $w\widetilde{w}^{-1}=_{_G} 1$,
% $$
%\area(w)=\area(w\widetilde{w}^{-1}).
% $$
%Since $T[e]$ is always the trivial path, this definition is compatible with the previous notion of area for closed
%paths at $e$. Finally, we define the \emph{area} of an arbitrary path $\gamma$ in $\Gamma(G)$ starting at $u$
%(non-necessarily equal to $e$) and ending at $v$, as the area of its translation $\gamma_e =u^{-1}\gamma$, which starts
%at $e$ and ends at $u^{-1}v$:
% $$
%\area(\gamma )=\area(\gamma_e) =\area(\gamma_e \widetilde{\gamma_e }^{-1}) =T[e,\, \iota \gamma]\,
%\gamma\widetilde{\gamma}^{\, -1} T[e,\, \iota \gamma]^{-1}.
% $$
%\note{QUINA DEF ES LA BONA ?????}
%%$\gamma T[\iota \gamma,\, \tau \gamma]^{-1}$, or of $T[e,\, \iota
%%\gamma]\gamma T[\iota \gamma,\, \tau \gamma]^{-1} T[e,\, \iota \gamma]^{-1}$.

\medskip

To analyze the mean Dehn function of a group $G$, we have to evaluate the sum of areas of all words in $A^{\ast}$
mapping to 1 in $G$, and having a given length. That is, the sum of areas of all paths in $\Gamma(G)$ of a given
length, and closed at $e$. To do this, we will do inductive arguments that force us to consider more general sums, like
the sum of areas of all paths in $\Gamma(G)$ starting at $e$ and of a given length (...and being closed or not). The
following notation will be useful in order to manipulate these sums.

For a given set of paths $P$ starting at $e$ (i.e. a given $P\subseteq A^*$) we denote by $\ar_P$ the sum of areas of
paths in $P$, $\ar_P =\sum_{\gamma \in P} \area (\gamma)$. Specially, if $v$ is a vertex in $\Gamma(G)$ and $n$ is a
positive integer, we denote by $\ar_v (n)$ the sum of areas of all paths $\gamma$ in $\Gamma(G)$ having length $n$,
starting at $\iota\gamma =e$ and ending at $\tau\gamma =v$. Note that, if $\lgr{v}>n$, then there are no such paths and
so $\ar_v (n)=0$. Note also that $\ar_{e}(n)$ is the sum of areas of all closed paths at $e$ with length $n$, which is
precisely the numerator of the spherical mean Dehn function of $G$ evaluated at $n$. Finally, let $\ar(n)$ denote the
sum of areas of all paths $\gamma$ in $\Gamma(G)$ having length $n$ and starting at $e$. Thus, we have
 $$
\begin{array}{rl}
\ar_v(n) & = \underset{\underset{\iota \gamma =e,\, \tau \gamma =v}{_{\la{\gamma}=n}}}{\sum} \area (\gamma),
\\ & \\
\ar(n) & =\underset{v}{\sum} \, \ar_v(n) = \underset{\underset{\iota \gamma =e}{_{\la{\gamma}=n}}}{\sum} \area
(\gamma),
\end{array}
 $$
Similarly, we denote by $\n_v(n)$ the number of paths $\gamma$ in $\Gamma(G)$ having length $n$, starting at
$\iota\gamma =e$ and ending at $\tau\gamma =v$. Of course, $\n_v (n)=0$ if $\lgr{v}>n$. Also, $\sum_v \n_v (n)=(2r)^n$.
This notation allows us to write
 $$
D_{\text{\rm smean}}(n)=\frac{\ar_e(n)}{\n_e(n)},
 $$
and suggests to define the \textit{open (spherical) mean Dehn function} as the averaged area over all such paths:
 $$
D_{\text{\rm osmean}}(n)=\frac{\ar(n)}{(2r)^n} =\frac{\sum_v \ar_v (n)}{\sum_v \n_v (n)}.
 $$

\medskip

In order to find an upper bound for $D_{\text{\rm osmean}}(n)$, we shall be guided by the following intuitive idea. Out
of the $(2r)^n$ paths of length $n$, those arriving ``far" from $e$ will mostly contribute with a ``big" area; but
there are ``few" of them. And those arriving ``close" to $e$ (which are ``much more" frequent) are going to contribute
less because they mostly have ``small" area.

To develop this intuitive idea, giving precise sense to the quoted words, we consider the following technical
condition. For all those finite presentations satisfying it, we will be able to give a recurrent estimation of
$\ar(n)$.

\begin{ass}\label{assumption} \emph{
Let $A=\{ a_1,\, \ldots ,\, a_r\}$, $F$ be the free group on $A$, and $G=\langle A\,|\, R\rangle$ be a finite
presentation of a quotient of $F$. For the rest of the present section we shall assume the existence of two
non-decreasing functions $f,g \colon \mathbb{N}\to \mathbb{R}^+$ and a constant $c_0$ such that, for every $c\gg 0$,
 $$
\sharp \{ w\in A^{\ast}\,|\, \la{w} =n,\,\, \lgr{w}>cf(n)\} =O\big( \frac{(2r)^n}{g(n)^{c-c_0}}\big).
 $$
(Note that this assumption is vacuous if $f(n)$ grows faster than linear, or if $c\leqslant c_0$.) }
\end{ass}

\begin{prop}\label{n-n/2}
Let $A=\{ a_1,\, \ldots ,\, a_r\}$, $F$ be the free group on $A$, and $G=\langle A\,|\, R\rangle$ be a finite
presentation of a quotient of $F$ satisfying assumption~\ref{assumption}. Choose an arbitrary geodesic combing $T$ in
$\Gamma(G)$. Then, for every $c\gg 0$, we have
 $$
\begin{array}{ll}
\ar(n)\leqslant & (2r)^{\lceil n/2 \rceil} \ar(\lfloor \frac{n}{2}\rfloor) + (2r)^{\lfloor n/2 \rfloor} \ar(\lceil
\frac{n}{2}\rceil) +\\ & \\ & (2r)^n D(4cf(n)) + D(2n) O\big( \frac{(2r)^n}{g(n)^{c-c_0}}\big).
\end{array}
 $$
\end{prop}

\demo Fix $c\in \mathbb{R}^+$ big enough from assumption~\ref{assumption}. Every summand in $\ar(n)$ has the form
$\area(\gamma )=\area(\gamma \widetilde{\gamma}^{-1})$ and so is bounded above by $D(2n)$ (since
$\la{\widetilde{\gamma}}\leqslant \la{\gamma}\leqslant n$). On the other hand, $\ar(n)$ is a sum of $(2r)^n$ summands.
Let us split $\ar(n)$ into two terms in such a way that we can improve one of these two estimates in each. Consider
$P_1 =\{ \gamma \,|\, \la{\gamma}=n,\, \iota(\gamma)=e,\, \lgr{\tau \gamma}>cf(n)\}$, $P_2 =\{ \gamma \,|\,
\la{\gamma}=n,\, \iota(\gamma)=e,\, \lgr{\tau \gamma}\leqslant cf(n)\}$. Separating
 \begin{equation}\label{1}
\ar(n) =\ar_{P_1}+ \ar_{P_2},
 \end{equation}
the first term has a small number of summands (according to assumption~\ref{assumption}), while the summands in the
second term are small (because they are areas of paths near to closed at $e$). More precisely,
 \begin{equation}\label{2}
\ar_{P_1}\leqslant D(2n)\cdot \sharp P_1 =D(2n)O\big( \frac{(2r)^n}{g(n)^{c-c_0}}\big),
 \end{equation}
and let us evaluate now the second term in~(\ref{1}). A typical summand there is the area of a path $\gamma$ of length
$n$, starting at $e$, and ending at some vertex $v$ such that $\lgr{v}\leqslant cf(n)$. That is, $\area (\gamma
\widetilde{\gamma}^{-1})$, where $\la{\gamma}=n$ and $\la{\widetilde{\gamma}}\leqslant cf(n)$. Break $\gamma$ into two
parts, $\gamma =\gamma_1 \gamma_2$ with $\la{\gamma_1 }=\lfloor \frac{n}{2}\rfloor$ and $\la{\gamma_2 }=\lceil
\frac{n}{2} \rceil$, and denote by $u$ the \emph{middle} point, $\tau \gamma_1 =u=\iota \gamma_2$ (see Figure~\ref{fig
1}, where $\widetilde{\gamma}_1=T[e,\, u]$, $\widetilde{\gamma}_2 =T[u,\, v]$ and $\widetilde{\gamma}=T[e,\, v]$).
\begin{figure}[htb]
%\input{1.pic}
%TeXCAD Picture [1.pic]. Options:
%\grade{\on}
%\emlines{\off}
%\epic{\off}
%\beziermacro{\on}
%\reduce{\on}
%\snapping{\off}
%\quality{8.00}
%\graddiff{0.01}
%\snapasp{1}
%\zoom{4.0000}
\unitlength 1mm % = 2.85pt
\linethickness{0.4pt}
\ifx\plotpoint\undefined\newsavebox{\plotpoint}\fi % GNUPLOT compatibility
\begin{picture}(88.5,30.75)(-2,65)
\qbezier(87.5,85.75)(60.63,90.75)(59.25,67.75) \put(67,83.5){\circle*{1}} \put(59.25,67.75){\circle*{1}}
\put(87.25,85.75){\circle*{1}}
%\emline(59,68.25)(59.25,68)
\multiput(59,68.25)(.03125,-.03125){8}{\line(0,-1){.03125}}
%\end
\put(59.25,68){\line(0,1){0}} \put(59.25,68){\line(1,2){7.75}} \put(67,83.5){\line(0,1){0}}
%\emline(67,83.5)(87,85.75)
\multiput(67,83.5)(.2985075,.0335821){67}{\line(1,0){.2985075}}
%\end
\put(87,85.75){\line(0,1){0}}
%\emline(87,85.75)(59.5,68.25)
\multiput(87,85.75)(-.052986513,-.03371869){519}{\line(-1,0){.052986513}}
%\end
\put(59.5,68.25){\line(0,1){0}}
%\emline(59.5,68.25)(59,68)
\multiput(59.5,68.25)(-.0625,-.03125){8}{\line(-1,0){.0625}}
%\end
\put(56.75,67.75){\makebox(0,0)[cc]{$e$}} \put(65.75,86){\makebox(0,0)[cc]{$u$}} \put(88.5,88){\makebox(0,0)[cc]{$v$}}
\put(60,79.25){\makebox(0,0)[cc]{$\gamma_1$}} \put(74.50,88.50){\makebox(0,0)[cc]{$\gamma_2$}}
\put(67,77){\makebox(0,0)[cc]{$\widetilde{\gamma}_1$}} \put(75.25,81.75){\makebox(0,0)[cc]{$\widetilde{\gamma}_2$}}
\put(74.5,75){\makebox(0,0)[cc]{$\widetilde{\gamma}$}}
%\vector(61.5,77.75)(62.75,80)
\put(63,80){\vector(1,2){.07}}
%\end
%\vector(75,86.5)(77.75,86.75)
\put(77.75,86.50){\vector(1,0){.07}}
%\end
%\vector(71.5,75.75)(73.75,77.5)
\put(73.75,77.25){\vector(4,3){.07}}
%\end
%\vector(63.5,76.5)(64.75,79.25)
\put(64.75,78.50){\vector(1,2){.07}}
%\end
%\vector(74.75,84.25)(76.75,84.75)
\put(76.75,84.75){\vector(4,1){.07}}
%\end
\end{picture}

    \caption{Breaking $\gamma$ into two parts.}
    \label{fig 1}
\end{figure}
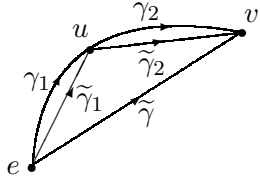
For every such $\gamma \in P_2$, we have
 $$
\begin{array}{rl}
\area(\gamma)=\area(\gamma_1 \gamma_2 \widetilde{\gamma}^{-1}) & \leqslant \area(\gamma_1 \widetilde{\gamma}_1^{\, -1})
+ \area(\widetilde{\gamma}_1 \gamma_2 \widetilde{\gamma}_2^{\,-1}
\widetilde{\gamma}_1^{\,-1}) + \area(\widetilde{\gamma}_1 \widetilde{\gamma}_2 \widetilde{\gamma}^{-1}) \\
& \\ & =\area(\gamma_1)+\area(\gamma_2) +\area(\widetilde{\gamma}_1 \widetilde{\gamma}_2 \widetilde{\gamma}^{-1}).
\end{array}
 $$
So,
 \begin{equation}\label{3}
\ar_{P_2} =\underset{\gamma\in P_2}{\sum} \area(\gamma) \leqslant \underset{\gamma\in P_2}{\sum}
(\area(\gamma_1)+\area(\gamma_2)) + \underset{\gamma\in P_2}{\sum} \area(\widetilde{\gamma}_1 \widetilde{\gamma}_2
\widetilde{\gamma}^{-1}).
 \end{equation}
To estimate the first summand in~(\ref{3}) observe that, moving $\gamma$ arround $P_2$, $\gamma_1$ moves inside the set
of words in $A^{\ast}$ of length $\lfloor \frac{n}{2} \rfloor$ (and $\gamma_2$ inside the set of words of length
$\lceil \frac{n}{2} \rceil$). Note also that every word of length $\lfloor \frac{n}{2} \rfloor$ appears as $\gamma_1$
at most $(2r)^{\lceil \frac{n}{2} \rceil}$ times (while every word of length $\lceil \frac{n}{2} \rceil$ appears as
$\gamma_2$ at most $(2r)^{\lfloor \frac{n}{2} \rfloor}$ times). Thus,
 \begin{equation}\label{4}
\underset{\gamma \in P_2}{\sum} (\area(\gamma_1)+\area(\gamma_2)) \leqslant (2r)^{\lceil \frac{n}{2} \rceil}
\ar(\lfloor \frac{n}{2} \rfloor) + (2r)^{\lfloor \frac{n}{2} \rfloor} \ar(\lceil \frac{n}{2} \rceil).
 \end{equation}
It remains to estimate the second summand in~(\ref{3}), i.e. the areas of geodesic triangles. To do this, we split
again $P_2$ into two disjoint sets, depending on $\lgr{u}$. Let $P_3 =\{ \gamma\in P_2 \,|\, \la{\widetilde{\gamma}_1
}=\lgr{u}>cf(\lfloor \frac{n}{2}\rfloor)\}$ and $P_4 =\{ \gamma\in P_2 \,|\, \la{\widetilde{\gamma}_1
}=\lgr{u}\leqslant cf(\lfloor \frac{n}{2}\rfloor)\}$, and
 \begin{equation}\label{5}
\underset{\gamma\in P_2}{\sum} \area(\widetilde{\gamma}_1 \widetilde{\gamma}_2 \widetilde{\gamma}^{-1}) =
\underset{\gamma\in P_3}{\sum} \area(\widetilde{\gamma}_1 \widetilde{\gamma}_2 \widetilde{\gamma}^{-1}) +
\underset{\gamma\in P_4}{\sum} \area(\widetilde{\gamma}_1 \widetilde{\gamma}_2 \widetilde{\gamma}^{-1}).
 \end{equation}
Again using the same argument as above, we can bound the first summand in~(\ref{5}) using the fact that it has few
summands,
 \begin{equation}\label{6}
\underset{\gamma\in P_3}{\sum} \area(\widetilde{\gamma}_1 \widetilde{\gamma}_2 \widetilde{\gamma}^{-1})\leqslant D(2n)
O\big( \frac{(2r)^n}{g(\lfloor \frac{n}{2}\rfloor)^{c-c_0}}\big).
 \end{equation}
Finally, the second summand in~(\ref{5}) can be bounded taking into account that all the involved triangles have
perimeter $\la{\widetilde{\gamma}_1}+\la{\widetilde{\gamma}_2}+\la{\widetilde{\gamma}}\leqslant
2(\la{\widetilde{\gamma}_1}+\la{\widetilde{\gamma}})\leqslant 2(cf(\lfloor \frac{n}{2}\rfloor)+cf(n))\leqslant 4cf(n)$.
Hence,
 \begin{equation}\label{7}
\underset{\gamma\in P_4}{\sum} \area(\widetilde{\gamma}_1 \widetilde{\gamma}_2 \widetilde{\gamma}^{-1})\leqslant
D(4cf(n))(2r)^n.
 \end{equation}
Combining together equations~(\ref{1}) to~(\ref{7}), we conclude the proof:
 $$
\begin{array}{rcl}
\ar (n) & \leqslant & (2r)^{\lceil \frac{n}{2} \rceil} \ar(\lfloor \frac{n}{2} \rfloor) + (2r)^{\lfloor \frac{n}{2}
\rfloor} \ar(\lceil \frac{n}{2} \rceil) + \vspace{.1cm} \\ & & D(2n)O\big( \frac{(2r)^n}{g(n)^{c-c_0}}\big) + D(2n)
O\big( \frac{(2r)^n}{g(\lfloor \frac{n}{2}\rfloor)^{c-c_0}}\big) + D(4cf(n))(2r)^n \vspace{.2cm} \\ \vspace{.1cm} & = &
(2r)^{\lceil \frac{n}{2} \rceil} \ar(\lfloor \frac{n}{2} \rfloor) + (2r)^{\lfloor \frac{n}{2} \rfloor} \ar(\lceil
\frac{n}{2} \rceil) + \vspace{.1cm} \\ & & (2r)^n D(4cf(n)) + D(2n) O\big( \frac{(2r)^n}{g(\lfloor
\frac{n}{2}\rfloor)^{c-c_0}}\big). \quad \Box
\end{array}
 $$
\bigskip

Let us make now another assumption to clear out one of the terms in the previous formula

\begin{ass}\label{assumption2} \emph{
From now on, we shall also assume that our group has polynomial Dehn function, say $D(n)=O(n^k)$ for some $k\in
\mathbb{R}^+$, and that our function $g(n)$ additionally satisfies that $\frac{g(n)}{n^{\alpha}}$ is uniformly bounded
away from zero, for some $\alpha >0$. }
\end{ass}

\begin{prop}\label{cleaned}
Under assumptions~\ref{assumption} and~\ref{assumption2}, we have
 $$
\ar(n)\leqslant (2r)^{\lceil n/2 \rceil} \ar(\lfloor \frac{n}{2}\rfloor) + (2r)^{\lfloor n/2 \rfloor} \ar(\lceil
\frac{n}{2}\rceil) + (2r)^n O(f(n)^k).
 $$
\end{prop}

\demo In the actual conditions, and taking $c>c_0 +k/\alpha$, the last term in the statement of Proposition~\ref{n-n/2}
will be
 $$
D(2n) O\big( \frac{(2r)^n}{g(n)^{c-c_0}}\big) \leqslant L\cdot (2r)^n \frac{n^k}{g(n)^{c-c_0}} \big(
\frac{g(n)}{n^{\alpha}} \big)^{c-c_0} =L(2r)^n n^{k+\alpha (c_0 -c)} \leqslant L (2r)^n,
$$
for an appropriate constant $L$, and so it is negligible:
 $$
\begin{array}{rcl}
\ar(n) & \leqslant & (2r)^{\lceil n/2 \rceil} \ar(\lfloor \frac{n}{2}\rfloor) + (2r)^{\lfloor n/2 \rfloor} \ar(\lceil
\frac{n}{2}\rceil) + (2r)^n D(4cf(n)) + D(2n) O\big( \frac{(2r)^n}{g(n)^{c-c_0}}\big) \\ & = & (2r)^{\lceil n/2 \rceil}
\ar(\lfloor \frac{n}{2}\rfloor) + (2r)^{\lfloor n/2 \rfloor} \ar(\lceil \frac{n}{2}\rceil) + (2r)^n O(f(n)^k). \quad
\Box
\end{array}
 $$
\bigskip

To conclude this section, let us unwrap the recurrence given at the previous statement, obtaining an upper bound for
the open spherical mean Dehn function of all finite presentations satisfying assumptions~\ref{assumption}
and~\ref{assumption2}.

\begin{thm}\label{open}
For every finite presentation (and geodesic combing) satisfying assumptions~\ref{assumption} and~\ref{assumption2}, and
for every non-decreasing function $h\colon \mathbb{N}\to \mathbb{R}^+$ satisfying $2h(\lceil \frac{n}{2}\rceil )+f(n)^k
\leqslant h(n)$ for $n\gg 0$, we have
 $$
D_{\text{osmean}}(n) =O\big( h(n)\big).
 $$
\end{thm}

\demo From Proposition~\ref{cleaned}, there exists a constant $M$ such that, for every $n\geqslant 2$,
 $$
\ar(n) \leqslant (2r)^{\lceil \frac{n}{2} \rceil} \ar(\lfloor \frac{n}{2} \rfloor) + (2r)^{\lfloor \frac{n}{2} \rfloor}
\ar(\lceil \frac{n}{2} \rceil) + M(2r)^n f(n)^k.
 $$
Now take $h$ as in the statement (for $n\geqslant n_0$), and let $K=\max \{M,\, \ar(2)/h(2),\ldots, \ar(n_0)/h(n_0)\}$.
Let us prove that, for $n\geqslant 2$,
 $$
\ar(n) \leqslant K(2r)^n h(n).
 $$
For $n=2,\ldots ,n_0$ the inequality is true, by construction. Fix a value of $n>n_0$, and assume the inequality true
for all smaller values. We have
 $$
\begin{array}{rcl}
\ar(n) & \leqslant & (2r)^{\lceil \frac{n}{2} \rceil} \ar(\lfloor \frac{n}{2} \rfloor) + (2r)^{\lfloor \frac{n}{2}
\rfloor} \ar(\lceil \frac{n}{2} \rceil) + M(2r)^n f(n)^k \\ & \leqslant & (2r)^{\lceil \frac{n}{2} \rceil}
K(2r)^{\lfloor \frac{n}{2} \rfloor} h(\lfloor\frac{n}{2} \rfloor) + (2r)^{\lfloor \frac{n}{2} \rfloor} K(2r)^{\lceil
\frac{n}{2} \rceil} h(\lceil \frac{n}{2} \rceil )+ M(2r)^n f(n)^k \\ & \leqslant & K(2r)^n \big( h(\lfloor\frac{n}{2}
\rfloor )+h(\lceil \frac{n}{2} \rceil ) +f(n)^k \big) \\ & \leqslant & K(2r)^n h(n).
\end{array}
 $$
Hence, $D_{\text{osmean}}(n) =\ar(n)/(2r)^n =O(h(n))$ concluding the proof. \qed

From Theorem~\ref{open} to being able to bound the spherical mean Dehn function, we will need to extract and use
another piece of information from the presentation of $G$. Namely, which proportion of the total $(2r)^n$ paths of
length $n$ are closed. Or, more generally, how sensible $\n_v(n)$ is in terms of $v$. This information strongly depends
on the group $G$ and on the specific presentation considered.

Finally, going from an estimation of the spherical mean Dehn function to an estimation of the mean Dehn function for
the same presentation, is easy after the following observation.

\begin{prop}\label{relation}
For any finite presentation of a group $G$, we have
 $$
D_{\text{mean}}(n) \leqslant  \underset{0\leqslant m\leqslant n}{\max} D_{\text{smean}}(m).
 $$
\end{prop}

\demo Directly from the definitions, we have
 $$
\underset {w\in B_G (n)}{\sum} \area(w)=\sum_{m=0}^n \underset {w\in S_G (m)}{\sum} \area(w) =\sum_{m=0}^n
D_{\text{smean}}(m) \cdot \sharp S_G (m) \leqslant
 $$
 $$
\leqslant \big( \underset{0\leqslant m\leqslant n}{\max} D_{\text{smean}}(m) \big) \sum_{m=0}^n \sharp S_G (m) =\big(
\underset{0\leqslant m\leqslant n}{\max} D_{\text{smean}}(m) \big) \cdot \sharp B_G (n). \quad \Box
 $$

\section{Counting words in abelian groups}\label{s-counting}

Let us apply now the techniques developed in the previous section to any finite presentation of an abelian group, until
obtaining explicit upper bounds for $D_{\text{osmean}}(n)$, $D_{\text{smean}}(n)$ and $D_{\text{mean}}(n)$. To do this,
we need first to verify that those presentations satisfy assumption~\ref{assumption} for appropriate functions $f,g$.
This is the goal of the present section.

\medskip

We start with a simple and well known lemma, which is straightforward to verify by induction.

\begin{lem}\label{fractions}
Let $x_1,\, \ldots ,\, x_r$ and $y_1,\, \ldots ,\, y_r$ be two lists of $r$ positive real numbers. Then,
 $$
\min \Big\{ \frac{x_1}{y_1},\, \ldots ,\, \frac{x_r}{y_r}\Big\} \leqslant \frac{x_1+\cdots +x_r}{y_1+\cdots + y_r
}\leqslant \max \Big\{ \frac{x_1}{y_1},\, \ldots ,\, \frac{x_r}{y_r}\Big\}. \Box
 $$
\end{lem}

Our arguments will strongly use the following lemma due to Kolmogorov (see Lemma 8.1 in page 378 of~\cite{M}). It seems
that this useful result proved in 1929, is somewhat forgotten in the literature and not known to many authors. For this
reason, and also for completion of the present paper, we add here a self-contained proof extracted from~\cite{M}. It
uses the following Tchebyshev inequality, which is straightforward to verify.

\begin{lem}[Tchebyshev]\label{Tch}
Let $X$ be a random variable and $f(x)$ be a nondecreasing real function. Then, for any real number $a$ such that
$f(a)> 0$, the following inequality holds:
 $$
{\text {\rm Pr}}\,(X>a)\leqslant \frac{E(f(X))}{f(a)}.
 $$
\end{lem}

\begin{lem}[Kolmogorov]\label{Kol}
Consider $n$ pairwise independent random variables $\{X_i\}$, $i=1, \ldots ,n$, with zero means and variances
$\sigma_i^2=E(X_i^2)$, and suppose that $|X_i| \leqslant d< \infty$. Let $S_n=\sum_{i=1}^nX_i$, and let $t$ be a real
number such that $0<td\leqslant s_n$, where $s_n^2=\var(S_n)=\sum_{i=1}^n\sigma_i^2$. Then, for any $\epsilon
>0$,
 $$
{\text{\rm Pr}}\,(S_n> \epsilon s_n)\leqslant \exp \bigl(-t\epsilon + \frac{1}{2}t^2 (1+\frac{1}{2}tds_n^{-1})\bigr).
 $$
\end{lem}

\demo For each $X_i$, and for every $j\geqslant 2$ we have
 $$
E(X_i^j)=E(X_i^{j-2}X_i^2 )\leqslant d^{\, j-2}E(X_i^2) =d^{\, j-2}\sigma_i^2.
 $$
Also, the following series are absolutely convergent and, since $0<tds_n^{-1}\leqslant 1$, and $\sum_{j=3}^{\infty}
\frac{2}{j!}=2(e-2.5)<0.5$, we have
 $$
\begin{array}{rl}
E\big( e^{ts_n^{-1}X_i}) & = E(\sum_{j=0}^{\infty} \frac{1}{j!}(ts_n^{-1}X_i )^j \big) \vspace{.2cm}\\ &
=\sum_{j=0}^{\infty } \frac{1}{j!}(ts_n^{-1})^j E(X_i^j) \vspace{.2cm}\\ & \leqslant 1+0+\sum_{j=2}^{\infty }
\frac{1}{j!} t^j s_n^{-j} d^{j-2}\sigma_i^2 \vspace{.2cm}\\ & =1+\frac{1}{2} (t\sigma_i s_n^{-1})^2 \big(
\sum_{j=2}^{\infty } \frac{2}{j!} (tds_n^{-1})^{j-2} \big) \vspace{.2cm}\\ & \leqslant 1+\frac{1}{2} (t\sigma_i
s_n^{-1})^2 \big( 1+tds_n^{-1} \sum_{j=3}^{\infty } \frac{2}{j!} (tds_n^{-1})^{j-3} \big) \vspace{.2cm}\\ & \leqslant
1+ \frac{1}{2}(t\sigma_i s_n^{-1})^2 (1+tds_n^{-1}\sum_{j=3}^{\infty} \frac{2}{j!}) \vspace{.2cm}\\ & \leqslant 1+
\frac{1}{2}(t\sigma_i s_n^{-1})^2 (1+\frac{1}{2}tds_n^{-1}) \vspace{.2cm}\\ & \leqslant \exp \big(
\frac{1}{2}(t\sigma_i s_n^{-1})^2 (1+\frac{1}{2}tds_n^{-1}) \big).
\end{array}
 $$
Now, using Tchebyshev's inequality (Lemma~\ref{Tch}) applied to $X=S_n$, $f(x)=e^{ts_n^{-1}x}$ and $a=\epsilon s_n$, we
have
 $$
\begin{array}{rl}
{\text {\rm Pr}}\,\big( S_n >\epsilon s_n \big) & \leqslant e^{-t\epsilon} E\big( e^{ts_n^{-1}S_n }\big) \vspace{.1cm} \\
& =e^{-t\epsilon } E\big( \prod_{i=1}^n e^{ts_n^{-1}X_i }\big) \vspace{.1cm} \\ & =e^{-t\epsilon} \prod_{i=1}^n E\big(
e^{ts_n^{-1}X_i }\big) \vspace{.1cm} \\ & \leqslant e^{-t\epsilon}\prod_{i=1}^n \exp \big( \frac{1}{2}(t\sigma_i
s_n^{-1})^2 (1+\frac{1}{2}tds_n^{-1}) \big) \vspace{.1cm} \\ & =\exp \big( -t\epsilon + \sum_{i=1}^n
\frac{1}{2}(t\sigma_i s_n^{-1})^2 (1+\frac{1}{2}tds_n^{-1}) \big) \vspace{.1cm} \\ & =\exp \big( -t\epsilon +
\frac{1}{2}t^2 (1+\frac{1}{2}tds_n^{-1}) \big).
\end{array}
 $$
This completes the proof. \qed

As a corollary, we easily deduce the following result on 1-dimensional random walks.

\begin{prop}\label{pabbigc}
Let $A=\{ a\}$ and let $F=G\simeq \mathbb{Z}$ be the infinite cyclic group generated by $A$. Given a real number $c>0$,
the number of words $w\in A^{\ast}$ with $\la{w}=n$ and $\lz{w}>c\sqrt{n\ln n}$ is $O(\frac{2^n}{n^{c-\frac12}})$.
\end{prop}

\demo Let us assume $n\geqslant 2$, and consider a 1-dimensional random walk on $\mathbb{Z}$ of length $n$, i.e. $n$
independent (and uniform) random variables $\{X_i\}$ with $X_i\in \{-1,1\}$ and $E(X_i)=0$, $i=1, \ldots ,n$. We have
$\sigma_i^2=1$ and $s_n^2=n$. Now, apply Kolmogorov Lemma with $d=1$, $t=\sqrt{\ln n}$ and $\epsilon=c\sqrt{\ln n}$. We
obtain that
 $$
\begin{array}{rcl}
{\text {\rm Pr}}\,\big( \sum_{i=1}^n X_i > c\sqrt{n\ln n}\big) & \leqslant & \exp \bigl( -c\ln n+\frac{\ln
n}{2}(1+\frac{\sqrt{\ln n}}{2\sqrt{n}} )\bigr) \\ & = & \exp \bigl( (\ln n) (-c+\frac{1}{2}+\frac{1}{4}\sqrt{\frac{\ln
n}{n}}\,) \bigr)
\\ & = & \dfrac{n^{\frac{1}{4}\sqrt{\frac{\ln n}{n}}}}{n^{c-\frac{1}{2}}} \vspace{.2cm} \\ & \leqslant &
\dfrac{K}{n^{c-\frac{1}{2}}},
\end{array}
 $$
where the last inequality is due to the fact that $\lim_{n\to \infty} n^{\frac{1}{4}\sqrt{\frac{\ln n}{n}}} =1$ (we can
take, for example, $K=1.35$).

But the number of words in $A^{\ast}$ of $A$-length $n$ is $2^n$. So, the previous inequality means that the number of
words $w\in A^{\ast}$ with $\la{w}=n$, $\lz{w}>c\sqrt{n\ln n}$, and representing positive integers is less than or
equal to $K\frac{2^n}{n^{c-\frac12}}$. By symmetry, the number of words $w\in A^{\ast}$ with $\la{w}=n$ and
$\lz{w}>c\sqrt{n\ln n}$ is at most $2K\frac{2^n}{n^{c-\frac12}}$. Finally, since $K$ does not depend on $n$ (neither on
$c$) we have the result. \qed

\medskip

The next statement is the analog of Proposition~\ref{pabbigc} for an arbitrary finitely generated abelian group.

\begin{prop}\label{pabbig}
Let $A=\{ a_1,\, \ldots ,\, a_r\}$, $F$ be the free group on $A$, and $G=\langle A\,|\, R\rangle$ be a finite
presentation of an abelian quotient of $F$. Given a real number $c>1/2$, the number of words $w\in A^{\ast}$ with
$\la{w}=n$ and $\lgr{w}>rc\sqrt{n\ln n}$ is $O\big( \frac{(2r)^n}{(\sqrt{n\ln n}\,)^{c-\frac12}}\big)$.
\end{prop}

\demo Since $G$ is an $r$-generated abelian group, the map $F\twoheadrightarrow G$ factors through $\mathbb{Z}^r$, so
we have $A^{\ast}\twoheadrightarrow F\twoheadrightarrow\mathbb{Z}^r \twoheadrightarrow G$. And, as we have observed
before, $\lzr{w} \geqslant \lgr{w}$. Therefore, it is enough to prove the result for $\mathbb{Z}^r$. So, we are reduced
to consider only the case where $G$ is the free abelian group of rank $r$.

%Observe that $rc\sqrt{n\ln n}\leqslant n-1$ for big enough $n$. Hence, we may assume that $n$ satisfies $rc\sqrt{n\ln
%n}\leqslant n-1$.
%
%WHY DO WE NEED THIS ??
%

Let $w\in A^{\ast}$. For any $i=1,\, \ldots ,\, r$, let $w_{a_i} \in \{ a_i\}^{\ast}$ be the word which can be obtained
from $w$ by deleting all letters different from $a_i$ and $a_i^{-1}$. Clearly, $\la{w}=\sum_{i=1}^r \la{w_{a_i}}$ (note
that $\la{w_{a_i}}=|w_{a_i}|_{\{a_i\}}$). Also, since $G$ is free abelian, $\lzr{w}= \sum_{i=1}^r \lzr{w_{a_i}}$.

Now, let $\ell=c\sqrt{n\ln n}$. %which we are assuming is less than or equal to $n$ WHY DO YOU NEED THIS ?????.
Note also that $\lzr{w}>r\ell$ implies $\lzr{w_{a_i}}>\ell$ for some $i$. Therefore, we have
 $$
\frac{\sharp \{w\in S(n) \mid \lzr{w}> r\ell\, \}}{\sharp S(n)}\leqslant \frac{\sum_{i=1}^r \sharp \{w\in S(n) \mid
\lzr{w_{a_i}} >\ell\,\}}{\sharp S(n)}.
 $$
Furthermore, for every $i=1,\, \ldots ,\, r$, we also have
 $$
\begin{array}{rcl}
\dfrac{\sharp \{ w\in S(n) \mid \lzr{w_{a_i}}>\ell\}}{\sharp S (n)} & = & \dfrac{\sum_{m=\lceil \ell\rceil }^{n}
\sharp \{w\in S(n) \mid \la{w_{a_i}} =m ,\, \lzr{w_{a_i}} \geqslant \lceil \ell\rceil\}}{\sharp S(n)} \vspace{.2cm}\\
& \leqslant & \dfrac{\sum_{m=\lceil \ell\rceil }^{n} \sharp \{w\in S(n) \mid \la{w_{a_i}} =m ,\, \lzr{w_{a_i}}
\geqslant \lceil \ell\rceil \}}{\sum_{m=\lceil \ell\rceil }^{n} \sharp \{w\in S(n) \mid \la{w_{a_i}}=m\}}
\vspace{.2cm}\\ & \leqslant & \underset{\lceil \ell\rceil \leqslant m\leqslant n}{\max} \dfrac{\sharp \{w\in S(n) \mid
\la{w_{a_i}} =m ,\, \lzr{w_{a_i}} \geqslant \lceil \ell\rceil \}}{\sharp \{w\in S(n) \mid \la{w_{a_i}}=m\}},
\end{array}
 $$
where the last inequality is justified by Lemma~\ref{fractions}. But, given a word $v\in \{ a_i\}^{\ast}$, the number
of words $w\in S(n)$ such that $w_{a_i}=v$ do not depend on $v$, but only on $m=\la{v}=|v|_{_{\{a_i\}}}$. So, for every
$\lceil \ell\rceil \leqslant m\leqslant n$, we have
 $$
\begin{array}{rcl}
\dfrac{\sharp \{w\in S(n) \mid \la{w_{a_i}} =m ,\, \lzr{w_{a_i}} \geqslant \lceil \ell\rceil \}}{\sharp \{w\in S(n)
\mid \la{w_{a_i}}=m\}} & = & \dfrac{\sharp \{v\in \{a_i\}^{\ast} \mid |v|_{_{\{a_i\}}}=m ,\,\, \lz{v}\geqslant \lceil
\ell\rceil \}}{\sharp \{v\in \{a_i\}^{\ast} \mid |v|_{_{\{a_i\}}}=m\}} \\ & = & \dfrac{\sharp \{v\in \{a_i\}^{\ast}
\mid |v|_{_{\{a_i\}}}=m ,\,\, \lz{v}> c\sqrt{n\ln n}\, \}}{2^m} \\ & \leqslant & \dfrac{\sharp \{v\in \{a_i\}^{\ast}
\mid |v|_{_{\{a_i\}}}=m ,\,\, \lz{v}> c\sqrt{m\ln m}\, \}}{2^m} \\ & \leqslant & \frac{K}{m^{c-\frac{1}{2}}},
\end{array}
 $$
for an appropriate constant $K$ (according to Proposition~\ref{pabbigc}, we can take $K=2.7$). Thus, collecting all
together,
 $$
\frac{\sharp \{w\in S(n) \mid \lzr{w}>rc\sqrt{n\ln n}\, \}}{(2r)^n}=\frac{\sharp \{w\in S(n) \mid \lzr{w}>r\ell\,
\}}{\sharp S(n)}
 $$
 $$
\leqslant \sum_{i=1}^r \Big( \underset{\lceil \ell\rceil \leqslant m\leqslant n}{\max} \, \frac{K}{m^{c-\frac{1}{2}}}
\Big) = \frac{rK}{\lceil \ell\rceil^{c-\frac{1}{2}}} \leqslant \frac{rK}{(c\sqrt{n\ln n}\, )^{c-\frac{1}{2}}},
 $$
where we used $c>1/2$. This proves that the number of words $w\in A^{\ast}$ with $\la{w}=n$ and $\lgr{w}>rc\sqrt{n\ln
n}$ is $O\big( \frac{(2r)^n}{(\sqrt{n\ln n}\,)^{c-\frac12}}\big)$.
%(a valid constant being $\frac{2.7r}{c^{c-\frac12}}$).
\qed

We can rephrase Proposition~\ref{pabbig} by saying that finite presentations of abelian groups satisfy
assumption~\ref{assumption}.

\begin{cor}
Let $A=\{ a_1,\, \ldots ,\, a_r\}$, $F$ be the free group on $A$, and $G=\langle A\,|\, R\rangle$ be a finite
presentation of an abelian quotient of $F$. The functions $f(n)=(n\ln n)^{1/2}$ and $g(n)=(n\ln n)^{1/2r}$ and the
constant $c_0=r/2$ satisfy assumption~\ref{assumption} for all $c>r/2$.
\end{cor}

\demo For any given $c>r/2$, Proposition~\ref{pabbig} tells us that
 $$
\begin{array}{rcl}
\sharp \{ w\in A^{\ast}\,|\, \la{w} =n,\,\, \lgr{w}>cf(n)\} & = & \sharp \{ w\in A^{\ast}\,|\, \la{w} =n,\,\,
\lgr{w}>r(c/r)\sqrt{n\ln n}\} \\ & = & O\big( \frac{(2r)^n}{(\sqrt{n\ln n}\,)^{\frac{c}{r}-\frac12}}\big) \\ & = &
O\big( \frac{(2r)^n}{g(n)^{c-\frac r2}}\big) \\ & = & O\big( \frac{(2r)^n}{g(n)^{c-c_0}}\big).
\end{array}
 $$
Hence, assumption~\ref{assumption} is satisfied starting at $c>r/2$. \qed

\section{The mean Dehn function of abelian groups}\label{s-ab-case}

The next step is to fulfill assumption~\ref{assumption2} for finite presentations of abelian groups. This is easy since
it is well known that those groups have quadratic Dehn function (take, $k=2$ in~\ref{assumption2}) and because
$g(n)=(n\ln n)^{1/2r}$ so, taking $\alpha =1/2r$, we have $\frac{g(n)}{n^{\alpha}}$ uniformly bounded away from zero.

In this situation, Theorem~\ref{open} allows us to deduce the following upper bound for the open spherical mean Dehn
function of an abelian group.

\begin{thm}\label{open-ab}
Let $A=\{ a_1,\, \ldots ,\, a_r\}$, $F$ be the free group on $A$, and $G=\langle A\,|\, R\rangle$ be a finite
presentation of an abelian quotient of $F$. Then,
 $$
D_{\text{osmean}}(n) =O\big( n(\ln n)^2 \big).
 $$
\end{thm}

\demo In our situation, Theorem~\ref{open} ensures us that $D_{\text{osmean}}(n) =O\big( h(n) \big)$ for every
non-increasing function $h\colon \mathbb{N}\to \mathbb{R}^+$ satisfying $2h(\lceil \frac{n}{2}\rceil )+n\ln n \leqslant
h(n)$ for $n\gg 0$. And this is the case of the function $h(n)=n(\ln n)^2$. An straightforward calculus exercise shows
that
 $$
2\lceil \frac{n}{2}\rceil (\ln \lceil \frac{n}{2}\rceil )^2 +n\ln n \leqslant 2\,\frac{n+1}{2}(\ln \frac{n+1}{2})^2
+n\ln n \leqslant n(\ln n)^2
 $$
is true, precisely for $n\geqslant 15$ (in fact, one can show that any function growing asymptotically more slowly does
not satisfy the required inequality). \qed

As announced at the end of Section~\ref{s-combings}, to estimate the spherical mean Dehn function, we need some more
information from the presentation of $G$, namely how the terms $\n_v(n)$ depend on the vertex $v$. For abelian groups,
this can be deduced from the following more general result.

\begin{thm} {\rm \cite[Chapter VI.5]{Var}}\label{nilp}.
Let $A=\{ a_1,\, \ldots ,\, a_r\}$, $F$ be the free group on $A$, and $G=\langle A\,|\, R\rangle$ be a finite
presentation of a virtually nilpotent quotient of $F$. Then,
 $$
\max_{v\in \Gamma(G)} \{ \n_v (n) \}=O\big( \frac{(2r)^n}{n^{d/2}}\big),
 $$
where $d$ is the degree of the (polynomial) growth function of $G$. Moreover, there exists another constant $L>0$ such
that
 $$
\n_e (n)\geqslant L\, \frac{(2r)^n}{n^{d/2}},
 $$
for every even $n\geqslant 2$.
\end{thm}

Regardless the meaning of $d$ (which is very significant within the group $G$ but is not relevant for the present
computations) the previous result allows us to transfer our upper bound to the spherical mean Dehn function.

\begin{thm}\label{spherical-ab}
Let $A=\{ a_1,\, \ldots ,\, a_r\}$, $F$ be the free group on $A$, and $G=\langle A\,|\, R\rangle$ be a finite
presentation of an abelian quotient of $F$. Then,
 $$
D_{\text{smean}}(n) =O\big( n(\ln n)^2 \big).
 $$
\end{thm}

\demo Using the present notation, we have $ D_{\text{\rm smean}}(n)=\frac{\ar_e(n)}{N_e(n)}$. We are going to estimane
the numerator again by cutting paths on two halfs. Let $P$ be the set of closed paths in $\Gamma(G)$, based at $e$ and
having length $n$. As in the proof of Proposition~\ref{n-n/2}, break every $\gamma \in P$ into two parts, $\gamma
=\gamma_1 \gamma_2$ with $\la{\gamma_1 }=\lfloor \frac{n}{2}\rfloor$ and $\la{\gamma_2 }=\lceil \frac{n}{2} \rceil$,
and denote by $u$ the \emph{middle} point, $\tau \gamma_1 =u=\iota \gamma_2$. We have
 $$
\area(\gamma) =\area(\gamma_1 \gamma_2 )\leqslant \area(\gamma_1 \widetilde{\gamma}_1^{\,-1}
)+\area(\widetilde{\gamma}_1 \gamma_2) = \area(\gamma_1) + \area(\gamma_2).
 $$
Now, taking into account that $\lgr{u} \leqslant \lfloor \frac n2 \rfloor$, and applying Theorem~\ref{nilp}, we have
 $$
\begin{array}{rl}
\ar_e(n) & =\underset{\underset{\phantom{a}}{\gamma \in P}}{\sum} \area(\gamma) \\ &
\leqslant \underset{\underset{\phantom{a}}{\gamma \in P}}{\sum} \area(\gamma_1) + \underset{\gamma \in P}{\sum} \area(\gamma_2 ) \\
&
= \underset{\underset{\phantom{a}}{0\leqslant \lgr{u}\leqslant \lfloor n/2 \rfloor}}{\sum} \ar_u (\lfloor \frac{n}{2}
\rfloor ) \n_u (\lceil \frac{n}{2} \rceil ) + \underset{0\leqslant \lgr{u}\leqslant \lfloor n/2 \rfloor}{\sum} \n_u
(\lfloor \frac{n}{2} \rfloor )\ar_u (\lceil \frac{n}{2} \rceil ) \\ &
\leqslant \underset{u\in \Gamma (G)}{\max} \{ \n_u (\lceil n/2\rceil) \} \cdot
\underset{\underset{\phantom{a}}{0\leqslant \lgr{u}\leqslant \lfloor n/2 \rfloor}}{\sum} \ar_{u}(\lfloor \frac{n}{2}
\rfloor )+ \underset{u\in \Gamma (G)}{\max} \{ \n_u (\lfloor n/2 \rfloor) \} \cdot
\underset{\underset{\phantom{a}}{0\leqslant \lgr{u}\leqslant \lfloor n/2 \rfloor}}{\sum} \ar_{u}(\lceil \frac{n}{2}
\rceil )
\\ &
= \underset{\underset{\phantom{a}}{u\in \Gamma (G)}}{\max} \{ \n_u (\lceil n/2\rceil) \} \cdot \ar (\lfloor \frac{n}{2}
\rfloor )+ \underset{u\in \Gamma (G)}{\max} \{ \n_u (\lfloor n/2 \rfloor ) \} \cdot \ar (\lceil \frac{n}{2} \rceil )
\\ &
\leqslant M \frac{(2r)^{\lceil n/2\rceil}}{\lceil n/2\rceil^{d/2}} \ar (\lfloor \frac{n}{2} \rfloor )+ M
\frac{(2r)^{\lfloor n/2 \rfloor}}{\lfloor n/2 \rfloor^{d/2}} \ar (\lceil \frac{n}{2} \rceil ),
\end{array}
 $$
for a appropriate constant $M$. Finally, applying again Theorem~\ref{nilp}, and Theorem~\ref{open-ab}, and collecting
together all the constants, we conclude
 $$
\begin{array}{rcl}
D_{\text{smean}}(n) =\displaystyle\dfrac{\ar_e(n)}{\n_e (n)} & \leqslant & \frac{M \frac{(2r)^{\lceil
n/2\rceil}}{\lceil n/2\rceil^{d/2}} \ar (\lfloor \frac{n}{2} \rfloor )+ M \frac{(2r)^{\lfloor n/2 \rfloor}}{\lfloor n/2
\rfloor^{d/2}} \ar (\lceil \frac{n}{2} \rceil )}{L\, \frac{(2r)^n}{ \underset{\phantom{a}}{n^{d/2}} }} \\ & \leqslant &
\frac{M}{L} \big( \frac{n}{\underset{\phantom{a}}{\lfloor \frac n2 \rfloor}} \big)^{d/2} \frac{(2r)^{\lceil n/2\rceil}
\ar (\lfloor \frac{n}{2} \rfloor )+ (2r)^{\lfloor n/2 \rfloor} \ar (\lceil \frac{n}{2} \rceil )}{(2r)^n} \\ & \leqslant
&
\frac{M}{L} \cdot 3^{d/2} \big( \frac{\ar (\lfloor \frac{n}{2} \rfloor )}{\underset{\phantom{a}}{(2r)^{\lfloor \frac n2
\rfloor}}} + \frac{\ar (\lceil \frac{n}{2} \rceil )}{(2r)^{\lceil \frac n2 \rceil }}\big)
\\ & = &
K \big( D_{\text{osmean}}(\underset{\phantom{a}}{\lfloor \frac n2 \rfloor })+ D_{\text{osmean}}(\lceil \frac{n}{2}
\rceil ) \big)
\\ & = &
O\big( n(\ln n)^2 \big).
\end{array}
 $$
However, a remark about the parity of the closed paths in $\Gamma(G)$ needs to be done here, since we have used the
second part of Theorem~\ref{nilp} for an arbitrary $n$, while it was stated only for the even ones. If all the
relations $R$ in our presentation have even length, then all closed paths have also even length, and
$D_{\text{smean}}(n)=0$ for every odd $n$, by convention. In this case, the above computations form a complete proof of
the Theorem, understanding everywhere that $n$ is even.

Otherwise, let $\gamma_0$ be a closed path in $\Gamma(G)$ of the smallest possible odd length, say $n_0$. Then for
every closed path $\gamma$ of even length $n$, $\gamma_0 \gamma$ is again a closed path, now of odd length $n+n_0$.
This proves that $\n_{e} (n+n_0) \geqslant \n_e (n)\geqslant L\, \frac{(2r)^n}{n^{d/2}}$. Adjusting the constants
appropriately, this shows that the assumption ``$n$ even" in the second part of Theorem~\ref{nilp} can be removed in
this case. Hence, the proof is complete. \qed

Finally, a similar result is true for the mean Dehn function.

\begin{thm}\label{mean-ab}
Let $A=\{ a_1,\, \ldots ,\, a_r\}$, $F$ be the free group on $A$, and $G=\langle A\,|\, R\rangle$ be a finite
presentation of an abelian quotient of $F$. Then,
 $$
D_{\text{mean}}(n) =O\big( n(\ln n)^2 \big).
 $$
\end{thm}

\demo This follows immediately from Theorem~\ref{spherical-ab} and Proposition~\ref{relation}, since $n(\ln n)^2$ is an
increasing function. \qed

\section*{Appendix}

At the end of the introduction, we pointed out the question of which is the most appropriate or natural notion of mean
Dehn function of a group $G=\langle A\,|\, R\rangle$, from the group theory point of view. That is, which is the set
that must be considered to average the areas over it? In this appendix we defend the opinion that the most appropriate
one is the set of closed paths in the Cayley graph $\Gamma(G,A)$ \textit{without backtrackings}, that is the set of
genuine words in the free group on $A$, mapping to 1 in $G$. However, we also want to illustrate that counting those
paths (and averaging the areas over them) seems to be a much more difficult task than doing the same over the set of
closed paths with possible backtrackings (as done in the present paper), or over the set of lazy words (as done
in~\cite{Y}).

Let $A$ be a finite set and $G=\langle A\,|\, R\rangle$ be a finite presentation of a group $G$. Above $G$ we can
consider the following tower of algebraic structures, each being a quotient of the previous one:
 $$
(A\cup \{1\})^{\ast}\twoheadrightarrow A^{\ast}\twoheadrightarrow F\twoheadrightarrow G.
 $$
Here, $(A\cup\{1\})^{\ast}$ is the free monoid on $A\cup A^{-1}\cup \{1\}$, $A^{\ast}$ is the free monoid on $A\cup
A^{-1}$, $F$ is the free group on $A$, and the arrows represent the canonical maps. The elements of these three
algebraic structures can be geometrically viewed into the Cayley graph $\Gamma(G,A)$: elements in $F$ (usually called
\textit{words}) are paths in $\Gamma(G,A)$ starting at $e$ and having no backtrackings; elements in $A^{\ast}$ (called
\textit{non-necessarily reduced words}) are paths in $\Gamma(G,A)$ starting at $e$ and having possible backtrackings;
finally, elements in $(A\cup\{1\})^{\ast}$ (called \textit{lazy words}) are paths in $\Gamma(G,A)$ starting at $e$,
with possible backtrackings, and allowed to temporarily stop at some of the visited vertices (one can think of them as
regular paths in the Cayley graph $\Gamma(G,A\cup \{1\})$, i.e. $\Gamma(G,A)$ with loops labeled 1 added everywhere).
Then, a path of each of these three types represents an element mapping to 1 in $G$ if and only if it is closed.

The intrinsic definition of area is for words mapping to 1 in $G$ (i.e. elements in the kernel of $F\twoheadrightarrow
G$). And the area of such a word is the minimal number of relations (again words) that are needed to express it. Then,
going to $F$ through the maps $(A\cup \{1\})^{\ast}\twoheadrightarrow A^{\ast}\twoheadrightarrow F$, the notion of area
naturally extends to non-necessarily reduced words, and to lazy words. Averaging then over length $n$ elements in these
three different sets, we get three different notions of mean Dehn function. From this point of view, the most natural
and canonical one seems to be that working directly in $F$, that is, averaging areas of words rather than
non-necessarily reduced or lazy words.

A completely different issue is the fact that averaging and estimating areas of words, even just counting words, seems
to be much more complicated and technically difficult that doing the same with non-necessarily reduced words, or with
lazy words. In this appendix we want to stress this difficulty by making some initial considerations about counting or
asymptotically estimating the number of closed paths without backtracking in the two dimensional integral lattice: a
timid and superficial starting into a field thats looks both interesting and complicated.

Before, we would like to suggest two more possible definitions of mean Dehn functions. If, for technical reasons, one
prefers to work with non-necessarily reduced words, then it makes sense to modify the notion of area by adding also the
number of cancelations needed. That is, think $G$ not as a quotient of the (free) group $F$ but as a quotient of the
(free) monoid $A^{\ast}$; then look at the monoid presentation $G=\langle a_1,\ldots ,a_r \,|\, R\cup \{ a_i a_i^{-1},
a_i^{-1}a_i \, | \, i=1,\ldots ,r \}\rangle$ and define, accordingly, the area of a word $w\in A^{\ast}$ with $w\eqg 1$
as the minimal number of relations in this monoid presentation required to express it. Averaging these areas over all
elements in $A^{\ast}$ of a given prefixed length, we get a new notion of mean Dehn function.

Similarly, we can also think $G$ as a quotient of the (free) monoid $(A\cup \{1\})^{\ast}$, then look at the monoid
presentation $G=\langle 1, a_1,\ldots ,a_r \,|\, R\cup \{1\}\cup \{ a_i a_i^{-1}, a_i^{-1}a_i \, | \, i=1,\ldots ,r
\}\rangle$ and define, accordingly, the area of a lazy word $w\in (A\cup \{1\})^{\ast}$ with $w\eqg 1$ as the minimal
number of relations in this monoid presentation required to express it (so additionally counting the number of 1's,
i.e. the total time lost in the corresponding random walk). Averaging these new areas over all elements in $(A\cup
\{1\})^{\ast}$ of a given prefixed length, we get another notion of mean Dehn function.

It seems interesting to analyze the relations between all these notions, and to understand up to which point they are
all equivalent, and independent of the presentation (if they are). We hope that future research works will clarify this
picure.

\bigskip

let $A=\{ a_1,\, \ldots ,\, a_r\}$ be an alphabet with $r$ letters, let $G=\langle A\,|\, R\rangle$ be a finite
presentation of a group $G$, and let $\Gamma=\Gamma(G,A)$ be the corresponding Cayley graph. Let $g_n$ be the number of
paths of length $n$ in $\Gamma$ which are closed at $e$ (denoted $\n_e(n)$ in section~\ref{s-combings}). And let $f_n$
be the total number of those having no backtracking. Clearly, $f_n\leqslant g_n$. Let us introduce generating functions
for $f_n$ and $g_n$:
 $$
F(t)=\sum_{n=0}^{\infty}f_nt^n,\quad \quad
G(t)=\sum_{n=0}^{\infty}g_nt^n.
 $$
The following formula connects $F(t)$ and $G(t)$ (see \cite{B}):
 \begin{equation}\label{bartholdi}
F(t)=\frac{1-t^2}{1+(2r-1)t^2}\cdot G\Bigl(\frac{t}{1+(2r-1)t^2}\Bigr).
 \end{equation}

\medskip

Let us concentrate now on the free abelian group of rank 2 with the standard set of $r=2$ generators, $G=\mathbb{Z}^2$
and $A=\{a,b\}$. And let us find both, exact formulas and the asymptotic behavior, for the corresponding numbers $f_n$.
It is clear that $f_n =g_n =0$ whenever $n$ is odd. So, we can restrict our attention to even lengths.

It is not difficult to see that $g_{2n}=\binom{2n}{n}^2$. Here is a very elegant argument that V. Guba pointed out to
us during his stay at CRM, Barcelona, in late 2004. A path of length $2n$ closed at the origin, is a sequence of $2n$
symbols from the alphabet $\{ a, a^{-1}, b, b^{-1} \}$ such that the total number of $a$'s coincide with that of
$a^{-1}$'s, and the total number of $b$'s coincide with that of $b^{-1}$'s. Consider the set of positions in the
sequence, $\{ 1, 2, \ldots ,2n \}$, and choose two subsets $C$ and $D$, both of cardinality $n$. Clearly, $\#
(C\setminus (C\cap D))=\# (D\setminus (C\cap D))$ and $\# (C\cap D) =\# (\{ 1, 2, \ldots ,2n \}\setminus (C\cup D))$.
We can then built a closed path at the origin by putting, for instance, $a$'s at the positions in $C\cap D$, $a^{-1}$'s
at the positions in $\{ 1, 2, \ldots ,2n \}\setminus (C\cup D)$, $b$'s at the positions in $C\setminus (C\cap D)$, and
$b^{-1}$'s at the positions in $D\setminus (C\cap D)$. This procedure gives a bijection between the set of paths we are
interested in, and the set of pairs of subsets $\{C,D \}$ of $\{ 1, 2, \ldots ,2n \}$ with cardinality $n$. Hence,
$g_{2n}=\binom{2n}{n}^2$.

So, retaking generating functions, we have $G(t)=\sum_{n=0}^{\infty} \binom{2n}{n}^2t^{2n}$. We shall use this to give
exact recurrent formulas for $f_{2n}$. Particularizing formula~(\ref{bartholdi}) to our case, we have
 $$
F(t)=\frac{1-t^2}{1+3t^2}\cdot G\Bigl(\frac{t}{1+3t^2}\Bigr).
 $$
Consider the following expansion
 $$
h(t)=\frac{t}{1+3t^2}=t-3t^3+9t^5-27t^7+\dots =\sum_{\stackrel{i=1}{i\,\,\text{\scriptsize odd}}}^{\infty}
(-3)^{\frac{i-1}{2}}t^i,
 $$
and denote by $A_{2n}$ the coefficient of $t^{2n}$ in the series
 $$
G(h(t))=1+\binom{2}{1}^2h(t)^2+\binom{4}{2}^2h(t)^4+\binom{6}{3}^2h(t)^6+\dots
 $$
Clearly $A_0=1$, and
 $$
A_{2n}=\binom{2}{1}^2A_{2n}^{(2)}+\binom{4}{2}^2A_{2n}^{(4)}+\binom{6}{3}^2A_{2n}^{(6)}+\dots +
\binom{2n}{n}^2A_{2n}^{(2n)},
 $$
where $A_{2n}^{(k)}$ is the coefficient of $t^{2n}$ it $h(t)^k$.

For any two natural numbers $m$ and $l$, denote by $P(m,l)$ the set of all ordered $l$-tuples $(i_1,i_2,\dots ,i_l)$
such that each $i_j$ is an odd positive number and $i_1+i_2+\dots +i_l=m$. We have
 $$
\begin{array}{rcl}
A_{2n}^{(2k)}t^{2n} & = & \displaystyle\sum_{(i_1,\dots ,i_{2k})\in P(2n,2k)} (-3)^{\frac{i_1-1}{2}}t^{i_1}\dots
(-3)^{\frac{i_{2k}-1}{2}}t^{i_{2k}} \\ & & \\ & = & \displaystyle\sum_{(i_1,\dots ,i_{2k})\in P(2n,2k)}
(-3)^{n-k}t^{2n}.
\end{array}
 $$
Now observe that, every $i_j$ in any $2k$-tuple from $P(2n,2k)$ is odd and so at least 1; hence, $\# P(2n,2k)$ equals
the number of ways of assigning the remaining $\frac{2n-2k}{2}=n-k$ (indistinguishable) twos into $2k$ boxes, namely
$\binom{n-k+2k-1}{2k-1}=\binom{n+k-1}{k+k-1}$. Thus
 $$
A_{2n}^{(2k)}=\# P(2n,2k) \cdot (-3)^{n-k} =\binom{n+k-1}{k+k-1}(-3)^{n-k},
 $$
and so,
 $$
A_{2n}=\sum_{k=1}^{n}\binom{2k}{k}^2 \binom{n+k-1}{k+k-1}(-3)^{n-k}.
 $$

Finally, since
 $$
 \frac{1-t^2}{1+3t^2}=1-4\sum_{s=1}^{\infty}(-3)^{s-1}t^{2s},
 $$
we have
 $$
f_{2n}=A_{2n}-4\sum_{s=1}^{n}(-3)^{s-1}A_{2n-2s}.
 $$
From this, we can deduce the following recurrent formula to compute the numbers $f_n$:
 $$
f_{2n}+3f_{2n-2}=A_{2n}-A_{2n-2}.
 $$

\medskip

For the problem of finding the asymptotic behavior of $f_n$ and $g_n$, define the numbers
 $$
\alpha=\underset{n\rightarrow \infty}{\overline {\lim}}\, f_n^{1/n},\quad \quad \beta =\underset{n\rightarrow
\infty}{\overline{\lim}}\,g_n^{1/n}
 $$
($\alpha$ is called the {\it co-growth} of the pair $(G,A)$, and $\frac{1}{2r}\beta$ the {\it spectral radius} of
$(G,A)$). In~\cite{Gr}, R.I. Grigorchuck found the following interesting formula relating $\alpha$, $\beta$ and $2r$
(the size of the alphabet):
 $$
\beta=\begin{cases} \alpha+\frac{2r-1}{\alpha} & {\text{if}\,\,\,\,\,} \alpha>\sqrt{2r-1},\\ & \\
\frac{2\sqrt{2r-1}}{2r} & {\text{otherwise}}.
\end{cases}
 $$
(Since $F(t)$ and $G(t)$ have radii of convergence $\frac{1}{\alpha}$ and $\frac{1}{\beta}$ respectively,
formula~(\ref{bartholdi}) connects the numbers $\alpha,\beta$ and $2r$, wherefrom one can deduce Grigorchuck's
formula.)

Back to the case of $\mathbb{Z}^2$, we have $g_{2n}=\binom{2n}{n}^2$. Hence, using Stirling's formula, $g_{2n}\sim
\frac{2}{\pi}\frac{4^{2n}}{2n}$. Thus, $\beta=4$. Since $2r=4$, Grigorchuk's formula implies that $\alpha=3$. Therefore
one can expect that $f_{2n}= O(\frac{3^{2n}}{2n})$. And using a result of Sharp, we prove that this is precisely the
asymptotic behavior of these numbers.

In~\cite{Sh}, R. Sharp gave an asymptotic formula for counting paths without backtrackings in the case $G=\mathbb{Z}^r$
and with respect to the standard set of generators. For $v\in G$, let $\mathcal{N}_v^{\, \prime}(n)$ denote the number
of paths in the Cayley graph, without backtrackings, having length $n$, starting at $e$ and ending at $v$. Consider
also the constant $\sigma$ given by
 $$
\sigma^2=\frac{1}{\sqrt{2r-1}}\Bigl[ 1+\Bigl( \frac{r+\sqrt{2r-1}}{r-\sqrt{2r-1}}\Bigr)^{1/2}\Bigr]
=\frac{\sqrt{2r-1}+1}{r-1}.
 $$

\begin{thm}[Sharp, \cite{Sh}]
Let $G\cong \mathbb{Z}^r$ be the free abelian group on $r\geqslant 2$ generators. With the above notations we have that
 $$
\lim_{\underset{n\in 2\mathbb{Z}}{n\rightarrow \infty }} \Bigl| \sigma^r n^{r/2}\cdot \frac{\mathcal{N}_v^{\,
\prime}(n)}{(2r)(2r-1)^{n-1}}-\frac{2}{(2\pi)^{r/2}}e^{-||v||^2/(2\sigma^2n)}\Bigr|=0,
 $$
uniformly in $v\in \mathbb{Z}^r$.
\end{thm}

The following two corollaries can be easily deduced from this result.

\begin{cor}
With the above notation for $\mathbb{Z}^2$,
 $$
f_{2n}\sim\frac{4}{3(\sqrt{3}+1)\pi}\cdot \frac{3^{2n}}{2n}.
 $$
\end{cor}

\medskip

\begin{cor}
There exist positive constants $C_1$ and $C_2$ (depending only on $r$) such that
 $$
C_1\cdot (2r-1)^{n} n^{-r/2}\leqslant \mathcal{N}_0^{\, \prime}(n)\leqslant \max_{v\in \mathbb{Z}^r}\,
\mathcal{N}_v^{\, \prime}(n)\leqslant C_2\cdot (2r-1)^{n} n^{-r/2}
 $$
for all positive even $n$.
\end{cor}

\medskip

The analysis performed above allowed us to obtain recurrent formulas for $f_{2n}$ and also its asymptotic behavior, in
the case of dimension 2. However, it is unclear to us how to use this information in order to obtain a good enough
estimate from above for the number of paths without backtracking, having length $n$, starting at $e$ and terminating
outside the ball of radius $\sqrt{n\ln n}$. Being able to do this, we would have the starting point to develop a
project similar to the one contained in the present paper, but centered on genuine words rather than non-necessarily
reduced words.

\section*{Acknowledgments}
We wish to thank V. Guba for fruitful discussions, especially concerning the appendix; in particular, the elegant
argument for counting $g_{2n}$ is due to him. Both authors thank the Centre de Recerca Matem\`atica at Barcelona for
the warm hospitality received during the fall semester of 2004, while most of this paper was done. The first named
author gratefully acknowledges support by the grant of the Complex integration projects of SBRAS N 1.9. and by the
INTAS grant N 03-51-3663. The second named author gratefully acknowledges partial support by DGI (Spain) through grant
BFM2003-06613.

\end{document}